\documentclass[11pt]{amsart}

\usepackage{amsmath}
\usepackage{amsxtra}
\usepackage{amscd}
\usepackage{amsthm}
\usepackage{amsfonts}
\usepackage{amssymb}
\usepackage{eucal}

\newcommand{\Sf}{\mathfrak{S}}
\newcommand{\vac}{{| {\rm vac} \rangle}}      
\newcommand{\sr}[1]{\scriptstyle{#1}}
\newcommand{\slth}{\widehat{\mathfrak{sl}}_2}
\newcommand\B{{\mc B}}
\newcommand{\bra}[1]{\langle #1 |}        
\newcommand{\ket}[1]{{| #1 \rangle}}      
\newcommand{\wt}{{\rm wt}\,}

\newcommand{\Hom}{\mathop{\rm Hom}}

\newcommand{\F}{\mathcal{F}}

\newcommand{\bea}{\begin{eqnarray}}
\newcommand{\ena}{\end{eqnarray}}
\newcommand{\be}{\begin{eqnarray*}}
\newcommand{\en}{\end{eqnarray*}}

\newcommand{\C}{{\mathbb C}}

\newcommand{\Z}{{\mathbb Z}}

\newcommand{\si}{{\sigma}}
\newcommand{\la}{{\lambda}}

\newcommand{\mc}{\mathcal}

\numberwithin{equation}{section}
\newtheorem{thm}{Theorem}[section]
\newtheorem{prop}[thm]{Proposition}
\newtheorem{lem}[thm]{Lemma}
\newtheorem{cor}[thm]{Corollary}

\newtheorem{dfn}[thm]{Definition}

\newtheorem{conj}[thm]{Conjecture}
\begin{document}

\title[Rigged paths]
{Sets of rigged paths with Virasoro characters}

\author{B. Feigin, M. Jimbo, T. Miwa, E. Mukhin and Y. Takeyama}
\address{BF: Landau institute for Theoretical Physics,
Chernogolovka,
142432, Russia}\email{feigin@feigin.mccme.ru}
\address{MJ: Graduate School of Mathematical Sciences,
The University of Tokyo, Tokyo 153-8914, Japan}
\email{jimbomic@ms.u-tokyo.ac.jp}
\address{TM: Department of Mathematics,
Graduate School of Science,
Kyoto University, Kyoto 606-8502,
Japan}\email{tetsuji@math.kyoto-u.ac.jp}
\address{EM: Department of Mathematics,
Indiana University-Purdue University-Indianapolis,
402 N.Blackford St., LD 270,
Indianapolis, IN 46202}\email{mukhin@math.iupui.edu}
\address{YT: Institute of Mathematics,
Graduate School of Pure and Applied Sciences,
University of Tsukuba, Tsukuba, Ibaraki 305-8571,
Japan}\email{takeyama@math.tsukuba.ac.jp}

\date{\today}

\begin{abstract}
Let $\{M^{(p,p')}_{r,s}\}_{1\le r\le p-1,1\le s\le p'-1}$ be
the irreducible Virasoro modules in the $(p,p')$-minimal series.
In our previous paper, we have constructed a monomial basis of
$\oplus_{r=1}^{p-1}M^{(p,p')}_{r,s}$ in the case $1<p'/p<2$.
By `monomials' we mean vectors of the form
$\phi^{(r_L,r_{L-1})}_{-n_L}
\cdots\phi^{(r_1,r_{0})}_{-n_1}\ket{r_0,s}$,
where $\phi_{-n}^{(r',r)}:M^{(p,p')}_{r,s}\to M^{(p,p')}_{r',s}$  
are the
Fourier components of the $(2,1)$-primary field and
$\ket{r_0,s}$ is the highest weight vector of $M^{(p,p')}_{r_0,s}$.
In this article, we introduce for all $p<p'$ with $p\geq3$ and $s=1$
a subset of such monomials as a conjectural basis of
$\oplus_{r=1}^{p-1}M^{(p,p')}_{r,1}$. We prove that
the character of the combinatorial set labeling these monomials
coincides with the character of the corresponding
Virasoro module.
We also verify the conjecture in the case $p=3$.
\end{abstract}
\maketitle

\section{Introduction}\label{sec:Intro}

Consider a representation $M$ of a vertex operator algebra $\mc V$.
Let $\phi^{(i)}(z)\in \mc V$ be a collection of fields,
and let $\phi_n^{(i)}$ be the corresponding Fourier coefficients.
\medskip

{\bf Question.} Find a set $I$ of sequences of pairs $(i_j,n_j)$ and a
vector $v\in M$ such that vectors $\{\phi_{n_L}^{(i_L)}\dots
\phi_{n_1}^{(i_1)}\phi^{(i_0)}_{n_0}v\ |\ (i_j,n_j)\in I\}$ form a
basis of $M$.
\medskip

Vectors of such form are usually called ``monomials''.
The question of finding a monomial basis, i.e., a basis consisting of
monomials, is an important, well-known and old problem,
solved in many interesting non-trivial cases.
Examples include: integrable $\widehat{\mathfrak{sl}}_n$ modules in  
terms of
$e_{1i}(z)$ currents \cite{LP} and \cite{Pr};
``big'' and ``small'' coinvariants of $\widehat{\mathfrak{sl}}_2$  
integrable
modules in terms of the currents $e(z), f(z), h(z)$ \cite{FKLMM1} and
\cite{FKLMM2};
${(2,2n+1)}$-Virasoro minimal series in terms of the Virasoro current
\cite{FF}; and $(3,3n\pm1)$-Virasoro minimal modules tensored with
Fock spaces in terms of an abelian current \cite{FJM}.

Following the same philosophy, in  \cite{FJMMT} we constructed a basis
of the $(p,p')$-Virasoro minimal series representations
$M^{(p,p')}_{r,s}$ ($1\leq r\leq p-1$, $1\leq s\leq p'-1$) with  
$1<p'/p<2$
and $p\geq3$, in terms of the $(2,1)$-primary field $\phi(z)$.
The basis has the form
\begin{equation}\label{KITEI}
\phi^{(r_L,r_{L-1})}_{-n_L}\cdots
\phi^{(r_{1},r_0)}_{-n_1}|b(s),s\rangle
\end{equation}
where $r_0=b(s),r_L=r, r_i=r_{i+1}\pm1$,
$\phi^{(r,r')}_{-n}:M^{(p,p')}_{r',s}\rightarrow M^{(p,p')}_{r,s}$
are the Fourier coefficients of $\phi(z)$,
and $n_i\in\Delta_{r_{i},s}-\Delta_{r_{i-1},s}+\Z$.
For each fixed $s$, $1\leq b(s)\leq p-1$ is so chosen
that the conformal dimension $\Delta_{r,s}$ of the space $M^{(p,p')} 
_{r,s}$
attains the minimum at $r=b(s)$.
We have also the condition
\begin{equation}\label{2-cond}
n_{i+1}-n_{i}\geq w(r_{i+1},r_{i},r_{i-1})\quad(1\leq i\leq L-1)
\end{equation}
where $w(r'',r',r)$ $(r=r'\pm1,r'=r''\pm1)$
are certain rational numbers (see \eqref{W1}--\eqref{W4} below).
It was shown that the condition (\ref{2-cond}) is a consequence of
the quadratic relations in the algebra generated
by Fourier coefficients of the primary field $\phi(z)$ .

In this paper we consider the case $p'/p>2$. In this case, in  
addition to
the quadratic relations, the algebra has cubic relations. To our  
surprise,
we found that in many cases (if not all) the cubic relations
result in a very simple exclusion rule in addition to (\ref{2-cond}):
\bea\label{3-cond}
n_{i+2}-n_i\ge 1.
\ena
The condition \eqref{3-cond} is void in the case $1<p'/p<2$
since it follows from (\ref{2-cond}).
We conjecture that in the general case $p<p'$
the monomials satisfying \eqref{2-cond} and \eqref{3-cond}
form a basis (see Conjecture \ref{CONJ} for the precise statement).

The purpose of this paper is to give some evidences for this  
conjecture.
Our results are two-fold. First, we show that the combinatorial set  
defined
by conditions (\ref{2-cond}), (\ref{3-cond}) (which we call the set of
`rigged paths') has the same graded character as that  of the  
Virasoro module
for any $(p,p')$ minimal theory with $s=1$
(see Theorem \ref{thm:main-theorem}). Second, we prove the conjecture
for the $(3,3n\pm1)$ minimal theories with $s=1$ (see Theorem \ref 
{thm:p=3}).

The first statement is purely combinatorial.
We prove it by showing that the characters of
Virasoro modules and the characters of the combinatorial sets
enjoy the same recurrence when the parameter changes from $p'-p$ to  
$p'$.
(see Proposition \ref{prop:t_to_t-1} and Proposition \ref 
{prop:ch_t_to_t-1}).

The second statement is based on the work \cite{FJM}, in which the
representation space
\be
V^{(3,p')}
=\Bigl(M^{(3,p')}_{1,1}\otimes(\oplus_{n:{\rm even}}\F_{n\beta})\Bigr)
\oplus\Bigl(M^{(3,p')}_{2,1}\otimes(\oplus_{n:{\rm odd}}\F_{n\beta}) 
\Bigr)
\en
of the abelian current $a(z)=\phi(z)\otimes \Phi_\beta(z)$ is studied.
Here $\Phi_\beta(z):\F_{n\beta}\rightarrow\F_{(n+1)\beta}$
is a certain bosonic vertex operator acting on the bosonic Fock space
$\F_{n\beta}$. It was shown that the abelian current $a(z)$ satisfies
cubic relations, and by exploiting the cubic relations
a monomial basis of the space $W^{(3,p')}$ generated from the vector
$|1,1\rangle\otimes |0\rangle$, where $|0\rangle\in\F_0$ is the  
highest weight
vector, was constructed. From the construction for $W^{(3,p'-3)}$,
we deduce a monomial basis of the space
\be
M^{(3,p')}=M^{(3,p')}_{1,1}\oplus M^{(3,p')}_{2,1}.
\en
The story is as follows.
We construct a filtration of the space $M^{(3,p')}$ by using the  
operator
$\phi(z)$. This filtration induces a current $\tilde a(z)$ acting on
the corresponding graded space, from the operator $a(z)$ acting on
$V^{(3,p')}$. The correlation functions for the operator $\tilde\phi 
(z)$
are equal to those for the operator $\tilde a(z)$ up to simple  
factors.
The latter belong to the space of correlation functions
of the operator $a(z)$ with $p'$ replaced by $p'-3$ up to simple  
factors.
The spanning property of the monomials \eqref{KITEI} with \eqref{2-cond},
\eqref{3-cond} is deduced by using these identifications of  
correlation
functions and the monomial basis of the space $W^{(3,p'-3)}$  
constructed
in\cite{FJM}.

The plan of our paper is as follows.
In Section \ref{sec:Rigged-path} we define the combinatorial set of  
rigged
paths and formulate Conjecture \ref{CONJ}.
In Section \ref{sec:Recurrence} we show that Virasoro characters
satisfy a recurrence relation when the parameter changes from $p'-p$
to $p'$. In Section \ref{sec:Thebijection} we show that the  
combinatorial
sets of rigged paths satisfy the same recurrence relation.
Section \ref{sec:p=3} is devoted to the proof of the conjecture
in the case $(p,p')=(3,3n\pm1)$.

In the paper by P. Jacob and P. Mathieu \cite{JM},
the problem of constructing monomial basis is studied,
and combinatorial conditions similar to ours are proposed.
Their study is restricted to the $(3,p)$ case, and in this case
the paths $(r_L,\ldots,r_1)$ do not appear in combinatorics.
We thank one of the referees for attracting our attentions to
this paper.

We also thank another referee for providing us with the simple proof
of Proposition \ref{REFEREE} as given in the below.
\setcounter{section}{1}
\setcounter{equation}{0}

\section{Rigged paths}\label{sec:Rigged-path}

\subsection{Minimal series}
We recall some definitions about minimal conformal field theory. For
the details we refer to \cite{DiFMS}.
Let {\it Vir} be the Virasoro algebra with the standard $\C$-basis
$\{L_n\}_{n\in\Z}$ and $c$ satisfying
\be
[L_m,L_n]=(m-n)L_{m+n}+\frac{c}{12}m(m^2-1)\delta_{m+n,0}\,,
\quad [c,L_n]=0\,.
\en
Fix a pair $(p,p')$ of relatively prime positive integers.
We set
\be
t=\frac{p'}{p}\,.
\en
We assume $p\geq3$ so that the $(2,1)$ primary field exists
(see below). We consider the minimal series of representations
$M^{(p,p')}_{r,s}$ ($1\leq r\leq p-1,1\leq s\leq p'-1$)
of {\it Vir}. The module $M^{(p,p')}_{r,s}$ is generated
by a vector $|r,s\rangle$ called the highest weight vector. It is an
irreducible module. The central element $c$ acts as the scalar
\be
c_{p,p'}=13-6\left(t+\frac{1}{t}\right)\,.
\en
The highest weight vector satisfies
\be
L_n|r,s\rangle=0\hbox{ if }n>0,
\quad L_0|r,s\rangle=\Delta^{(t)}_{r,s}|r,s\rangle\,,
\en
where
\[
\Delta^{(t)}_{r,s}=\frac{(rt-s)^2-(t-1)^2}{4t}\,.
\]

We fix $s$. The $(2,1)$ primary field
\be
\phi^{(r\pm1,r)}(z)=\sum_{n\in\Z-\Delta^{(t)}_{r,s}+\Delta^{(t)}_{r 
\pm1,s}}
\phi^{(r\pm1,r)}_{-n}z^{n-\Delta^{(t)}_{2,1}}
\en
is a generating series of linear operators $\phi^{(r\pm1,r)}_{-n}$  
acting as
\be
\phi^{(r\pm1,r)}_{-n}:M^{(p,p')}_{r,s}\rightarrow M^{(p,p')}_{r 
\pm1,s}\,.
\en
Up to a scalar multiple, they are characterized by the commutation  
relations
with the Virasoro generators:
\be
&&[L_n,\phi^{(r\pm 1,r)}(z)]=z^n\bigl(z\partial+(n+1)\Delta^{(t)}_ 
{2,1}\bigr)
\phi^{(r\pm 1,r)}(z)\,.
\en
The module $M^{(p,p')}_{r,s}$ is graded by eigenvalues of the  
operator $L_0$:
\be
M^{(p,p')}_{r,s}&=&\oplus_d\left(M^{(p,p')}_{r,s}\right)_d\,,
\\
(M^{(p,p')}_{r,s})_d&=&\{|v\rangle
\in M^{(p,p')}_{r,s} \mid
L_0|v\rangle=d|v\rangle\}\,.
\en
The character $\chi^{(p,p')}_{r,s}(q)$ is defined by the formula
\be
\chi^{(p,p')}_{r,s}(q)={\rm tr}_{M^{(p,p')}_{r,s}}q^{L_0}\,.
\en
The  primary field preserves the grading:
\be
&&\phi^{(r',r)}_{-n}(M^{(p,p')}_{r,s})_d\subset (M^{(p,p')}_{r',s})_ 
{d+n}\,.
\en

\subsection{Quadratic relations}
In \cite{FJMMT} we derived a set of quadratic relations for
the Fourier components of the $(2,1)$ primary field
of the following form:
\begin{equation}
\sum_{r'\atop n+n'=d}c^{(r,r',r'')}_{n,n'}\phi^{(r,r')}_{-n}
\phi^{(r',r'')}_{-n'}=0\,,
\label{QR}
\end{equation}
where $r,r''$ and $d$ are fixed, and the coefficients
$c^{(r,r',r'')}_{n,n'}$ are such that if $n'>N$ for some $N$
they are zero except for finitely many of them.
Let us describe a consequence of the quadratic relations
in a little bit different way than in \cite{FJMMT}.
Define weights $w^{(t)}(r,r',r'')=w^{(t)}(r'',r',r)=w^{(t)}(p-r,p- 
r',p-r'')$
where $r-r',r'-r''\in\{-1,1\}$ by
\begin{eqnarray}
w^{(t)}(r,r\pm1,r\pm2)&=&\frac t2,\quad\hbox{if }1\leq r,r\pm1,r\pm2 
\leq p-1,
\label{W1}\\
w^{(t)}(r,r+1,r)&=&2-\frac t2+[rt]-rt,\quad\hbox{if }2\leq r,r+1 
\leq p-1,
\label{W2}\\
w^{(t)}(r,r-1,r)&=&1-\frac t2-[rt]+rt,\quad\hbox{if }1\leq r-1,r 
\leq p-2,
\label{W3}\\
w^{(t)}(1,2,1)&=&w^{(t)}(p-1,p-2,p-1)=3-\frac{3t}2.
\label{W4}
\end{eqnarray}
Here $[x]$ is the integer part of $x$.
\begin{prop}\label{SPAN}
Any monomial of the form
\begin{equation}
\phi^{(r_L,r_{L-1})}_{-n_L}\cdots\phi^{(r_1,r_0)}_{-n_1}
\label{MONOM}
\end{equation}
can be written as an infinite linear combination of monomials
satisfying
\begin{equation}
n_{i+1}-n_i-w^{(t)}(r_{i+1},r_i,r_{i-1})\geq0.
\label{WA}
\end{equation}
\end{prop}
\begin{proof}
Essentially, the proof is given in \cite{FJMMT}
(see Proposition 3.3 in that paper). Note, however, that
in \cite{FJMMT} the weights $w^{(t)}(1,2,1)$ and $w^{(t)} 
(p-1,p-2,p-1)$ are
given as special cases of \eqref{W2} and \eqref{W3}, respectively.
This is because
the range of $t$ was restricted to $1<t<2$, wherein the two  
expressions
coincide. In \cite{FJMMT}, the quadratic relations were derived  
without
the restriction on $t$, and the statement that any monomial of the  
form
\eqref{MONOM} can be rewritten into an infinite linear combination
of those which satisfy \eqref{WA} was, in effect, proved without
the restriction $t<2$ by using \eqref{W4},
instead of using the special cases of \eqref{W2} and \eqref{W3}.
\end{proof}

We give a few more remarks.

As stated above, if we rewrite a monomial by using
the relations of the form \eqref{QR}, we obtain an infinite linear  
combination.
Note, however, that if
we fix the degree $n_1+\cdots+n_L$
and restrict $n_1\geq N$ for some $N$, there exists only finitely many
monomials satisfying the condition \eqref{WA}.
Therefore, such infinite sums are meaningful after completion.

Let us consider vectors of the form
\begin{equation}
\phi^{(r_L,r_{L-1})}_{-n_L}\cdots\phi^{(r_1,r_0)}_{-n_1}|r_0,s 
\rangle\,,
\label{VEC}
\end{equation}
where $|r_0,s\rangle$ is the highest weight vector of $M^{(p,p')}_ 
{r_0,s}$.
By abuse of terminology we call such vectors {\it monomials}.

Proposition \ref{SPAN} implies that any vector of the above form
belongs to the linear span (in the sense of finite linear  
combinations)
of those satisfying \eqref{WA} and the highest weight condition
\be
n_1+\Delta^{(t)}_{r_0,s}-\Delta^{(t)}_{r_1,s}\geq0.
\en

\subsection{Statement of the results}
Now we restrict to the case
\be
t>1 \hbox{ and }(r_0,s)=(1,1).
\en
We say a monomial
\begin{equation}
\phi^{(r_L,r_{L-1})}_{-n_L}\cdots\phi^{(r_1,r_0)}_{-n_1}|1,1\rangle
\label{VEC2}
\end{equation}
is {\it admissible} if and only if \eqref{WA}, and
\begin{equation}
n_1\geq\Delta^{(t)}_{2,1}=\frac{3t-2}4,
\label{HW}
\end{equation}
and
\begin{equation}
n_{i+2}-n_i\geq1
\label{AD}
\end{equation}
when $r_{i-1},r_i,r_{i+1}r_{i+2}\in\{r,r+1\}$
for some $1\leq r\leq p-2$, $1\leq i\leq L-2$, are satisfied.

Set
\begin{eqnarray}
v^{(t)}(r)&=&1-w^{(t)}(r,r+1,r)-w^{(t)}(r+1,r,r+1).
\label{VT}
\end{eqnarray}
We have
\begin{eqnarray}
v^{(t)}(r)&=&\begin{cases}p'-5&\hbox{if }p=3;\\[0pt]
[2t-3] &\hbox{if $p>3$ and $r=1,p-2$};\\[0pt]
[(r+1)t]-[rt]-2 &\hbox{if }1<r<p-2.
\end{cases}
\nonumber
\end{eqnarray}
Although $w^{(t)}$ is not necessarily an integer, $v^{(t)}$ are  
integers.
Note also that $v^{(t)}(r)=v^{(t)}(p-1-r)$.

{\it A rigged path of length} $L$ is a table of integers of the form
\be
P=
\begin{pmatrix}
r_L&r_{L-1}&\cdots&r_1&r_0\\
&\sigma_{L-1}&\cdots&\sigma_1&\sigma_0
\end{pmatrix}
\en
where $r_0=1$, $1\leq r_i\leq p-1$ $(0\leq i\leq L)$ and
\begin{equation}
r_{i+1}-r_i\in\{-1, 1\}\quad(0\leq i\leq L-1)\,.
\label{PATH}
\end{equation}
In the usual terminology, a sequence of integers $(r_i)$ satisfying
\eqref{PATH} is called a path. A rigged path is decorated by the  
rigging
$(\sigma_i)$. For brevity, we often call a rigged path simply a  
{\it path}.

A path is called admissible at level $t$ if
$\sigma_i\geq0\quad(0\leq i\leq L-1)$, and
\begin{equation}
\sigma_i+\sigma_{i+1}\geq v^{(t)}(r)
\label{AD2}
\end{equation}
when $r_{i-1},r_i,r_{i+1},r_{i+2}\in\{r,r+1\}$
for some $1\leq r\leq p-2$, $1\leq i\leq L-2$.

Note that these conditions correspond to
\eqref{HW}, \eqref{WA} and \eqref{AD}, respectively, if we set
\be
\sigma_0=n_1-\Delta_{2,1}^{(t)}, \quad
\sigma_i=n_{i+1}-n_i-w^{(t)}(r_{i+1},r_i,r_{i-1}) \quad
(1\le i \le  L-1).
\en

If $t<2$, we have $v^{(t)}(r)\leq0$ and the condition \eqref{AD2}  
follows
{}from the positivity of $\sigma_i$'s. If $t>2$, we have
\be
v^{(t)}(r)\geq0\quad(1\leq r\leq p-2),
\en
and moreover
\bea\label{VT1}
v^{(t)}(1)=v^{(t)}(p-2)\geq1.
\ena

We denote by $C^{(t)}_L$ the set of
rigged paths of length $L$ which are admissible at level $t$,
and by $C^{(t)}_{L,r}$ the subset of $C^{(t)}_L$ consisting of paths
such that $r_L=r$. The subset $C^{(t)}_{L,r}$ is empty unless
$r\equiv L+1\bmod2$.

We define the degree of $P\in C^{(t)}_L$ by
\bea
d(P)&=&\sum_{i=1}^Ln_i\label{DEGREE}\\
&=&L\Delta^{(t)}_{2,1}+\sum_{i=1}^{L-1}(L-i)w^{(t)}(r_{i+1},r_i,r_ 
{i-1})
+\sum_{i=0}^{L-1}(L-i)\sigma_i,\nonumber
\ena
and the character of $C^{(t)}_{L,r}$ by
\be
{\rm ch}_q\,C^{(t)}_{L,r}=\sum_{P\in C^{(t)}_{L,r}}q^{d(P)}.
\en
Our main result is the following identity:
\begin{thm}\label{thm:main-theorem}
We have
\be
\chi^{(p,p')}_{r,1}(q)=\sum_{L \ge 0}{\rm ch}_q\,C^{(t)}_{L,r}.
\en
\end{thm}
Theorem \ref{thm:main-theorem} will be proved
at the end of Section \ref{sec:Thebijection}.
This result motivates us to make the following conjecture:
\begin{conj}\label{CONJ}
For $1 \le r \le p-1$,
the set of admissible monomials of the form \eqref{VEC},
where $r_L=r, \, r_{0}=1$ and $s=1$, is a basis of
$M^{(p,p')}_{r,1}$.
\end{conj}
The case $1<t<2$ of the conjecture has been proved in \cite{FJMMT}.
We prove Conjecture \ref{CONJ} in the case $p=3$
in Section \ref{sec:p=3}.

\section{Recurrence structure}\label{sec:Recurrence}
\subsection{Recurrence relation for Virasoro characters}

Recall the formula for $\chi_{r,s}^{(p, p')}(q)$
\cite{RC}:
\begin{eqnarray}
\chi_{r,s}^{(p, p')}(q)=
\frac{q^{\Delta_{r,s}^{(t)}}}{(q)_{\infty}}\left(
\sum_{n \in \mathbb{Z}}q^{pp'n^{2}+(p'r-ps)n}-
\sum_{n \in \mathbb{Z}}q^{pp'n^{2}+(p'r+ps)n+rs}
\right).
\label{eq:bosonic-formula}
\end{eqnarray}
Here $(q)_{\infty}=\prod_{j=1}^{\infty}(1-q^{j})$.
In the case $1<t<2$,
$\chi_{r, 1}^{(p, p')}(q)$ was written
in the following form, \cite{FJMMT,W}:
\begin{eqnarray}
\chi_{r, 1}^{(p, p')}(q)=
q^{\Delta_{r, 1}^{(t)}}\sum_{m \equiv r-1 \, {\rm mod} \, 2}
\frac{1}{(q)_{m}}
K^{(p, p'-p)}_{m, r}(q),
\label{eq:fermionic-1<t<2}
\end{eqnarray}
where
\begin{eqnarray}
K_{m, r}^{(p, \bar{p'})}(q)=
q^{\frac{m^{2}-(r-1)^{2}}{4}}
\sum_{n \in \mathbb{Z}}
q^{p\bar{p'}n^{2}+\bar{p'}nr}
\left(
{m \atopwithdelims[] \frac{m-r+1}{2}-pn}-
{m \atopwithdelims[] \frac{m+r+1}{2}+pn} \right).
\label{eq:Kostka-part}
\end{eqnarray}
Here we have used the notation
\begin{eqnarray*}
(q)_{n}=\prod_{j=1}^{n}(1-q^{j}), \quad
{m \atopwithdelims[] n}=\frac{(q)_{m}}{(q)_{n}(q)_{m-n}}.
\end{eqnarray*}
The identity \eqref{eq:fermionic-1<t<2} can be generalized to the  
case $t>2$:
\begin{prop}\label{prop:recursion-character}
Set $k=[t]$. Then the following equality holds:
\begin{eqnarray}
\chi_{r,1}^{(p,p')}(q)=
q^{\Delta_{r,1}^{(t)}}
\sum_{m_{0}, \ldots , m_{k-1} \ge 0 \atop m_{0}\equiv r-1 \, {\rm  
mod} \,\, 2}
\frac{q^{Q^{(k)}(m_{0}, \ldots , m_{k-1})-\frac{k-1}{4}(r^{2}-1)}}
{(q)_{m_{0}} \cdots (q)_{m_{k-1}}}
K_{m_{0}, r}^{(p, p'-kp)}(q).
\label{eq:recursion-character}
\end{eqnarray}
Here
\begin{eqnarray*}
Q^{(k)}(m_{0}, \ldots , m_{k-1})&=&
\frac{k-1}{4}m_{0}^{2}+\sum_{j=1}^{k-1}(k-j)m_{j}^{2} \\
&+&
\sum_{j=1}^{k-1}(k-j)m_{0}m_{j}+
2\sum_{1 \le j<j' \le k-1}(k-j')m_{j}m_{j'} \\
&+&
\frac{k-1}{2}m_{0}+\sum_{j=1}^{k-1}(k-j)m_{j}.
\end{eqnarray*}
\end{prop}
Proposition \ref{prop:recursion-character} will be proved in
Section \ref{subsec:proof-of-recursion}.

For $L \in \mathbb{Z}_{\ge 0}$
denote by $\chi_{r, 1: L}^{(p, p')}(q)$
the right hand side of \eqref{eq:recursion-character} with
the sum replaced by the partial sum
over $m_{0}, \ldots , m_{k-1}$ satisfying $m_{0}+2(m_{1}+\cdots+m_ 
{k-1})=L$.
Then we have
\begin{eqnarray*}
\chi_{r, 1}^{(p, p')}(q)=\sum_{L \ge 0}\chi_{r, 1; L}^{(p, p')}(q).
\end{eqnarray*}
It follows from Proposition \ref{prop:recursion-character} that
\begin{prop}\label{prop:t_to_t-1}
We have
\begin{eqnarray}
\chi_{r, 1; L}^{(p, p')}(q)=
\sum_{m \ge 0}
\frac{q^{\frac{L^{2}}{4}+\frac{L}{2}}}{(q)_{m}}
\chi_{r, 1; L-2m}^{(p, p'-p)}(q).
\label{eq:recursion-str-1}
\end{eqnarray}
\end{prop}

In Section \ref{sec:Thebijection} we prove
(see Theorem \ref{thm:iota-is-bijection})
\begin{prop}\label{prop:ch_t_to_t-1}
\begin{eqnarray}
{\rm ch}_{q}C_{L, r}^{(t)}=
\sum_{m \ge 0}
\frac{q^{\frac{L^{2}}{4}+\frac{L}{2}}}{(q)_{m}}
{\rm ch}_{q}C_{L-2m, r}^{(t-1)}.
\label{eq:recursion-str-2}
\end{eqnarray}
\end{prop}
Theorem \ref{thm:main-theorem} follows from
these identities.
\subsection{Proof of Proposition \ref{prop:recursion-character}}
\label{subsec:proof-of-recursion}

In the following we set
$1/(q)_{n}=0$ for $n<0$.
To prove Proposition \ref{prop:recursion-character}
we use the following formula (the proof given here is provided by  
one of
the referees).:

\begin{prop}\label{REFEREE}
For $l \in \mathbb{Z}_{\ge 0}$ and $\mu \in \mathbb{Z}$, we have
\begin{eqnarray}
\frac{1}{(q)_{\infty}}=
\sum_{N_{0}, \ldots , N_{l} \ge 0}
\frac{q^{\sum_{j=0}^{l}N_{j}^{2}+\mu\sum_{j=0}^{l}N_{j}}}
{(q)_{N_{0}+\mu}(q)_{N_{0}}(q)_{N_{1}-N_{0}} \cdots (q)_{N_{l}-N_ 
{l-1}}}.
\label{eq:recursion-lemma-3}
\end{eqnarray}
\end{prop}

\begin{proof}
Consider the following functions $g_{l}(z) \,\, (l=0, 1, \ldots )$:
\begin{eqnarray*}
g_l(z)=\sum_{N_0,\dots,N_l \ge 0}\frac{z^{\sum_{j=0}^lN_j}q^{\sum_ 
{j=0}^lN_j^2}}
{(zq)_{N_0}(q)_{N_0}(q)_{N_1-N_0}\dots(q)_{N_l-N_{l-1}}},
\end{eqnarray*}
where $(z)_{n}:=\prod_{j=0}^{n-1}(1-zq^{j})$.
Let us prove that $g_{l}(z)=1/(zq)_{\infty}$.
Then we obtain \eqref{eq:recursion-lemma-3} by setting $z=q^{\mu}$.

First note that $g_{0}(z)$ can be rewritten in terms of
the basic hypergeometric series ${}_{2}\phi_{1}$ as follows:
\begin{eqnarray*}
g_0(z)=\lim_{a,b\rightarrow\infty}{}_2\phi_1(a, b ; zq ; q, zq/ab),
\end{eqnarray*}
where
\begin{eqnarray*}
{}_{2}\phi_{1}(a, b; c; q , z):=
\sum_{n=0}^{\infty}\frac{(a)_{n}(b)_{n}}{(c)_{n}(q)_{n}}\,z^{n}.
\end{eqnarray*}
On the other hand we have the $q$-Gauss sum identity (see, e.g.,  
\cite{GR}):
\begin{eqnarray*}
{}_{2}\phi_{1}(a, b; c; q, c/ab)=
\frac{(c/a)_{\infty}(c/b)_{\infty}}{(c)_{\infty}(c/ab)_{\infty}}.
\end{eqnarray*}
This implies that $g_{0}(z)=1/(zq)_{\infty}$.

Next consider the case where $l>0$.
Rewrite the definition of $g_{l}(z)$ as
\begin{eqnarray*}
g_{l}(z)=\sum_{N_1\dots,N_l \ge 0}
\frac{z^{\sum_{j=1}^lN_j}q^{\sum_{j=1}^lN_j^2}}
{(q)_{N_1}(q)_{N_1-N_0}\dots(q)_{N_l-N_{l-1}}}\,
f_{N_{1}}(z),
\end{eqnarray*}
where
\begin{eqnarray*}
f_{N}(z)=(q)_{N}\sum_{j=0}^{N}
\frac{z^{j}q^{j^{2}}}{(zq)_{j}(q)_{j}(q)_{N-j}}.
\end{eqnarray*}
It is easy to see that
\begin{eqnarray*}
f_{N}(z)=\lim_{a \to \infty}{}_{2}\phi_{1}(a, q^{-N}; zq; q, zq^{N 
+1}/a).
\end{eqnarray*}
Using the $q$-Gauss sum formula again, we see that
$f_{N}(z)=1/(zq)_{N}$.
This implies $g_{l}(z)=g_{l-1}(z)$ and
therefore we get $g_{l}(z)=g_{0}(z)=1/(zq)_{\infty}$.
\end{proof}

\medskip
\noindent{\it Proof of Proposition \ref{prop:recursion-character}}.

Substituting \eqref{eq:Kostka-part} to the right hand side of
\eqref{eq:recursion-character},
we obtain two sums, which we refer to as I and II.
First consider I:
\begin{eqnarray}
&&
\sum_{m_{0}, \ldots , m_{k-1} \ge 0 \atop m_{0}\equiv r-1 \, {\rm  
mod} \,\, 2}
\frac{q^{Q^{(k)}(m_{0}, \ldots , m_{k-1})-\frac{k-1}{4}(r^{2}-1)}}
{(q)_{m_{0}} \cdots (q)_{m_{k-1}}}
\label{eq:recursion-proof-1}
\\
&& \qquad {}\times
q^{\frac{m_{0}^{2}-(r-1)^{2}}{4}}
\sum_{n \in \mathbb{Z}}
q^{p(p'-kp)n^{2}+(p'-kp)nr}
{m_{0} \atopwithdelims[] \frac{m_{0}-r+1}{2}-pn}.
\nonumber
\end{eqnarray}
Set
\begin{eqnarray*}
N_{0}=\frac{m_{0}-r+1}{2}-pn, \quad
N_{j}-N_{j-1}=m_{j} \quad  (j=1, \ldots , k-1)
\end{eqnarray*}
and rewrite \eqref{eq:recursion-proof-1} as
a summation over $N_{0}, \ldots , N_{k-1} \ge 0$ and $n \in \mathbb 
{Z}$.
Then \eqref{eq:recursion-proof-1} becomes
\begin{eqnarray*}
\sum_{n \in \mathbb{Z}}q^{pp'n^{2}+(p'r-p)n}F_{k}(r-1+2pn),
\end{eqnarray*}
where
\begin{eqnarray*}
F_{k}(\mu)=
\sum_{N_{0}, \ldots , N_{k-1} \ge 0}
\frac{q^{\sum_{j=0}^{k-1}N_{j}^{2}+\mu N_{0}+
(\mu+1)\sum_{j=1}^{k-1}N_{j}}}
{(q)_{N_{0}+\mu}(q)_{N_{0}}(q)_{N_{1}-N_{0}} \cdots
(q)_{N_{k-1}-N_{k-2}}}.
\end{eqnarray*}

In the same way rewrite II
by setting
\begin{eqnarray*}
N_{0}=\frac{m_{0}+r+1}{2}+pn, \quad
N_{j}-N_{j-1}=m_{j} \quad (j=1, \ldots , k-1).
\end{eqnarray*}
The result is
\begin{eqnarray*}
\sum_{n \in \mathbb{Z}}q^{pp'n^{2}+(p'r+p)n+r}F_{k}(-r-1-2pn).
\end{eqnarray*}
Thus we find that the right hand side of \eqref{eq:recursion-character}
is equal to
\begin{eqnarray}
q^{\Delta_{r, 1}^{(t)}} \sum_{n \in \mathbb{Z}}
q^{pp'n^{2}+(p'r-p)n}\left(
F_{k}(r-1+2pn)-q^{r+2pn}F_{k}(-r-1-2pn) \right).
\label{eq:recursion-proof-2}
\end{eqnarray}
Set $\mu=r+2pn$. The last part of (\ref{eq:recursion-proof-2})
reads $F_{k}(\mu-1)-q^{\mu}F_{k}(-\mu-1)$.
In the sum $F_{k}(-\mu-1)$, we replace $N_j$ by $N_j+\mu$. Then, we  
have
\begin{eqnarray*}
&&
F_{k}(\mu-1)-q^{\mu}F_{k}(-\mu-1) \\
&& {}=
(1-q^{\mu})
\sum_{N_{0}, \ldots , N_{k-1} \ge 0}
\frac{q^{\sum_{j=0}^{k-1}N_{j}^{2}+\mu\sum_{j=0}^{k-1}N_{j}}}
{(q)_{N_{0}+\mu}(q)_{N_{0}}(q)_{N_{1}-N_{0}} \cdots
(q)_{N_{k-1}-N_{k-2}}}=
\frac{1-q^{\mu}}{(q)_{\infty}}.
\end{eqnarray*}
Here we used \eqref{eq:recursion-lemma-3} in the last equality.
Hence \eqref{eq:recursion-proof-2} is equal to
the right hand side of the formula \eqref{eq:bosonic-formula}.
This completes the proof.
\qed
\medskip

\section{The bijection}\label{sec:Thebijection}
\subsection{Particles}
In each path $P\in C^{(t)}_L$, we locate ``particles''.

Recall, that a path $P\in C^{(t)}_L$ is a table of integers
\bea
P=\left(\begin{matrix}
r_L& r_{L-1}&\dots&r_{1}&r_0\\
& \si_{L-1}&\dots&\si_1&\si_0
\end{matrix}\right).
\ena
The integers $r_x$  satisfy
the conditions $r_0=1$,
$1\leq r_x\leq p-1$ $(1\leq x \leq L)$ and
$|r_x-r_{x-1}|=1$ $(1\leq x \leq L)$. The integers
$\si_x$ satisfy the conditions
$\sigma_x\geq0$ $(0\leq x\leq L-1)$ and
\bea
\si_x+\si_{x-1}&\geq&v^{(t)}(r)\quad(2\leq x\leq L-1)
\label{B2}
\ena
if
\bea
&&r_{x+1},r_x,r_{x-1},r_{x-2}\in\{r,r+1\}\hbox{ for some }
1\leq r\leq p-2.\label{B3}
\ena
It is convenient to define $\sigma_L=\infty$, and use the convention
that $\infty\pm1=\infty$.

Our aim is to relate the sets $C_L^{(t-1)}$ ($L=0,1,2,\ldots$)
to the sets $C_L^{(t)}$ ($L=0,1,2,\ldots$).
There is an injective mapping from
$C_L^{(t-1)}$ to $C_L^{(t)}$ such that the image of the mapping
is equal to the subset of $C_L^{(t)}$ consisting of
paths $P$ for which the inequality (\ref{B2}) is strict,
and $\sigma_x\geq1$ ($1\leq x \leq L-1$), and $\sigma_x\geq2$
if $r_{x+1}=r_{x-1}=1$ or $p-1$. We say the path $P$ has no  
particle in such
case. When these conditions are violated,
we observe the appearance of ``particles'' in $P$.
We define the number of particles in $P$.
Then, we construct a bijection from $C_{L-2m}^{(t-1)}\times\pi_m$ to
the subset of paths in $C^{(t)}_L$ with $m$ particles.
Here we denote by $\pi_m$ the set of partitions $(\la_1,\ldots,\la_m)$
of length $m$, i.e.,
\be
\la_1\geq\la_2\geq\cdots\geq\la_m\geq0.
\en
We now give a precise definition of particles.

Given a path $P\in C^{(t)}_L$, we define an equivalence relation in  
the set
$\{1,\ldots,L-1\}$. We say a neighboring pair of integers
$x,x-1\in\{1,\ldots,L-1\}$ is connected if and only if \eqref{B3}  
is valid and
\bea\label{B4}
\si_x+\si_{x-1}&=&v^{(t)}(r)\quad(2\leq x\leq L-1).
\ena

We say $x\sim y$, where $x,y\in\{1,\ldots,L-1\}$,
if and only if all neighboring pairs of integers
in the interval between $x$ and $y$ are connected.

We call an equivalence class $B$ for this equivalence relation
{\it a block of particles} if one of the following is satisfied:
\bea\label{BLOCK}
&&|B|\geq2,\\
&&B=\{x\}\hbox{ and }\si_x=1, (r_{x+1},r_x,r_{x-1})=(1,2,1)\hbox{ or }
(p-1,p-2,p-1),\\
&&B=\{x\}\hbox{ and }\si_x=0, r_{x+1}=r_{x-1}.
\ena
We denote by $\B(P)$ the set of blocks for the path $P$.
We call a block $B$ an isolated particle if $B$ consists of one  
element,
$|B|=1$.

For each block $B\in\B(P)$,
let $\max(B),\min(B)$ be the largest and the smallest
integer in $B$.
The blocks in $P$ are naturally ordered:
$B>B'$ if and only if ${\rm min}(B)>{\rm max}(B')$.

If $x$ belongs to a block $B$, we have $r_{x+1}=r_{x-1}$.
We denote this number by $r'_x$.

\begin{dfn}\label{NUMPAR}
Let $P\in C^{(t)}_L$. We define a map $m:\B(P)\rightarrow\Z_{\geq0}$.
Namely, for $B\in\B(P)$, we set
\bea
m(B)=
\begin{cases}
\left[\frac{|B|}2\right]
&\hbox{ if }\sigma_x\geq2\hbox{ or if }
\sigma_x=1\hbox{ and }r'_x\not=1,p-1;\\[5pt]
\left[\frac{|B|+1}2\right]
&\hbox{ if }\sigma_x=0\hbox{ or if }\sigma_x=1\hbox{ and }r'_x=1,p-1.
\end{cases}
\ena
where $x=\max(B)$.
The number $m(B)$ is called the number of particles in the block $B$.
We also set
\bea
m(P)=\sum_{B\in\B(P)}m(B).
\ena
The number $m(P)$ is called the number of particles in the path $P$.
\end{dfn}
If $|B|$ is odd, then we have $\sigma_{{\rm min}(B)}=\sigma_{{\rm  
max}(B)}$.
Therefore, one can use $x={\rm min}(B)$ in Definition \ref{NUMPAR}.
Note also that if $m(P)\leq[L/2]$.
\subsection{Propagation of particles}
Roughly speaking, blocks in a path $P$ are the location of  
particles in
the sequence of integers
\be
L-1,L-2,\ldots,2,1.
\en
We move the particles in $P$ by changing blocks locally
in this sequence without changing the total number of particles.
We number the particles from the left to the right.
We define the left (resp., right) move of the $j$-th particle, $M^+_j$
(resp., $M^-_j$). These operations change a path $P\in C^{(t)}_L$ to
a path $M^\pm_jP\in C^{(t)}_L$. For some paths $P\in C^{(t)}_L$,  
the path
$M^\pm_jP\in C^{(t)}_L$ is not defined.
If a block contains more than one particles, i.e., $m(B)\geq2$,
the left move is defined only for the leftmost particle in the block,
and the right move is only for the rightmost one. In other words,
particles in a block are located without gaps and they cannot pass  
each other.
We start with the definition of the domain of the operations $M^ 
\pm_j$.

We set
\bea
D^+(P)&=&\{\max(B)|B\in\B(P)\},\\
D^-(P)&=&
\begin{cases}
\{\min(B)|B\in\B(P)\}\backslash\{1\}&\hbox{ if $\sigma_0=\sigma_1=0 
$};\\
\{\min(B)|B\in\B(P)\}&\hbox{ otherwise}.
\end{cases}
\ena
We define maps $m^\pm_P:D^\pm(P)\rightarrow\{1,\ldots,m(P)\}$:
\bea
m^+_P(x)=1+\sum_{B\in\B(P)\atop \max(B)>x}m(B),\\
m^-_P(x)=\sum_{B\in\B(P)\atop \min(B)\geq x}m(B).
\ena
Set
\bea
I^\pm(P)={\rm Im}(m^\pm_P)\subset\{1,\ldots,m(P)\}.
\ena
The map $m^\pm_P$ is injective, and therefore, on the image we can  
define
the inverse mappings
\bea
x^\pm_P=(m^\pm_P)^{-1}:I^\pm(P)&\rightarrow&\{1,\ldots,L-1\},\\
j&\mapsto&x_P^\pm(j).
\ena
The left move $M^+_jP$ is defined if and only if $j\in I^+(P)$, and
the right move $M^-_jP$ is defined if and only if $j\in I^-(P)$.
In such cases, we define the position of the $j$-th particle by $x^ 
+_P(j)$
or $x^-_P(j)$. If  $j\in I^+(P)\cap I^-(P)$, then we have
$x^+_P(j)=x^-_P(j)$. In this case, the corresponding particle is  
isolated.

In general, for each  $1\leq j\leq m(P)$ there exists a unique block
$B_j\in\B(P)$ such that
\bea
1+\sum_{B\in\B(P)\atop B>B_j}m(B)\leq j\leq
\sum_{B\in\B(P)\atop B\geq B_j}m(B).
\ena
We say that the $j$th particle is located in the block $B_j$.  
However, we
do not specify its position except for the leftmost one or the  
rightmost one.
\begin{dfn}
Let $P\in C^{(t)}_L$. Recall our convention $\sigma_L=\infty$.
Fix $j\in I^\pm(P)$ and set $x=x^\pm_P(j)$.
We define
\bea
M^\pm_jP
=\left(\begin{matrix}
r^\pm_L& r^\pm_{L-1}&\dots&r^\pm_{1}&r^\pm_0\\
& \si^\pm_{L-1}&\dots&\si^\pm_1&\si^\pm_0
\end{matrix}\right),
\ena
by setting $r^\pm_y=r_y$ $(0\leq y\leq L)$ and
$\si^\pm_y=\si_y$ $(0\leq y\leq L-1)$ except for the following.

{\it Case 1}\quad$\si_x\not=0:$
\bea
\si^\pm_x&=&\si_x-1,\\
\si^\pm_{x\mp1}&=&\si_{x\mp1}+1.
\ena

{\it Case 2}\quad$\si_x=0$ and $r'_x=1,p-2:$
\bea
\si^\pm_x&=&1,\\
\si^\pm_{x\pm1}&=&\si_{x\pm1}-1.
\ena

{\it Case 3}\quad$\si_x=0$ and $r'_x\not=1,p-2:$
\bea
r^\pm_x&=&2r'_x-r_x,\\
\si^\pm_{x\pm1}&=&
\begin{cases}
\si_{x\pm1}-v^{(t)}(r)-1&\hbox{ if }r_{x\pm2}=r_x;\\
\si_{x\pm1}+v^{(t)}(r+\varepsilon)&\hbox{ if }r_{x\pm2}=r_x+2 
\varepsilon,
\end{cases}\\
\si^\pm_{x\mp1}&=&
\begin{cases}
\si_{x\mp1}-v^{(t)}(r)&\hbox{ if }r_{x\mp2}=r_x;\\
\si_{x\mp1}+v^{(t)}(r+\varepsilon)+1&\hbox{ if }r_{x\mp2}=r_x+2 
\varepsilon,
\end{cases}
\ena
where $r={\rm min}\,\{r_x,r'_x\}$ and $\varepsilon=\pm1$.
\end{dfn}
The moves $M^\pm_i$ change the degree of a path by $\pm1$ and do  
not change
the number of particles:
\begin{lem}
Let $P\in C_{L,r}^{(t)}$ and $j\in I^\pm(P)$.
We have $M^\pm_jP\in C_{L,r}^{(t)}$,
$m(M^\pm_jP)=m(P)$ and $d(M^\pm_jP)=d(P)\pm1$.
\end{lem}
The proof is only case-checking.
\subsection{Properties of $M_j^\pm$}
We list properties of moves $M_j^\pm$.
In most cases, we omit proofs because they are only case-checkings.
However, the following remark might help understanding.

If $j\in I^+(P)$ and the $j$-th particle belongs to
a block $B$ with more than one particles, after the move $M^+_j$
this particle quits the block and the number of particles in $B$  
decreases.
The $j$-th particle either becomes an isolated particle or
joins in another block to increase the number of particles in that  
block.
The consideration of $M^-_j$ is similar.

The moves $M^+_j$ and $M^-_j$ are the inverse to each other:
\begin{lem}\label{INV}
Let $j\in I^\pm(P)$.
Then, we have $j\in I^\mp(M^\pm_jP)$ and $M^\mp_jM^\pm_jP=P$.
\end{lem}

The moves $M^\varepsilon_i$ and $M^\varepsilon_j$ are commutative
as far as they are defined:
\begin{lem}\label{++}
Suppose that $i\not=j+1$. If $M^\pm_iM^\pm_jP$ is defined, then
$M^\pm_jM^\pm_iP$ is also defined and $M^\pm_jM^\pm_iP=M^\pm_iM^ 
\pm_jP$.
\end{lem}

The moves $M^+_i$ and $M^-_j$ are commutative as far as they are  
defined:
\begin{lem}\label{+-=-+}
Suppose that $i\not=j$. If $M^\pm_iM^\mp_jP$ is defined, then
$M^\mp_jM^\pm_iP$ is also defined and $M^\mp_jM^\pm_iP=M^\pm_iM^ 
\mp_jP$.
\end{lem}

We are particularly interested in the relation between the moves
$M^+_{j+1}$ and $M^-_j$:
\begin{lem}\label{NEIGHBOR}
The move $M^+_{j+1}P$ is defined if and only if the move $M^-_jP$
is defined.

The move $M^+_{j+1}M^-_jP$ is defined if and only if the move
$M^-_jM^+_{j+1}P$ is defined, and in such a case we have
$M^+_{j+1}M^-_jP=M^-_jM^+_{j+1}P$.
\end{lem}

\begin{cor}\label{l+=l-}
 Let $P\in C_L^{(t)}$ and let $l$ be a positive integer.
The move $(M^+_{j+1})^lP$  is defined if and only if
the move $(M^-_j)^lP$ is defined.
\end{cor}
\begin{proof} We use induction on $l$. The case $l=1$ is proved in
Lemma \ref{NEIGHBOR}. Suppose that we have proved the statement for  
$l-1$.
By Lemma \ref{NEIGHBOR}
$(M_{j}^-)^lP=M_j^-(M^-_j)^{l-1}P$ is defined if and only if
$M^+_{j+1}(M^-_j)^{l-1}P$ is defined. Again by Lemma \ref{NEIGHBOR}
$M^+_{j+1}(M^-_j)^{l-1}P$ is defined if and only if
$M^{l-1}_jM^+_{j+1}P$ is defined. By induction hypothesis
$(M^-_j)^{l-1}M^+_{j+1}P$ is defined if and only if
$(M_{j+1}^+)^{l-1}M^+_{j+1} P$ is defined. Thus, we have proved that
$(M_{j+1}^+)^l P$
is defined if and only if the path $(M_{j}^-)^lP$ is defined.
\end{proof}
\subsection{Rigging of particles}
We define the rigging of particles $\lambda_j(P)$ $(j=1,\ldots,m(P))$.
Let $\pi_m$ be the set of partitions of length $m$.
The rigging $(\la_1(P),\la_2(P),\dots,\la_{m(P)}(P))$
is by definition an element of $\pi_{m(P)}$.

Let $P\in C_L^{(t)}$. We set formally $\lambda_{m(P)+1}(P)=0$.
Starting from $j=m(P)$ we define $\lambda_j(P)$ inductively
by requiring that $\la_j(P)-\la_{j+1}(P)$ is equal to
the maximal number $l$ such that the path $(M^-_j)^lP$ is defined.
Thus we obtain a partition.
We call $\la_j(P)$ {\it the rigging of the $j$-th particle} in the  
path $P$.

Note that by Corollary \ref{l+=l-} if $j\not=m(P)$
the number $\la_j(P)-\la_{j+1}(P)$ is also equal to the maximal  
number $l$
such that path $(M^+_{j+1})^lP$ is defined. Therefore, using Lemma  
\ref{INV},
we have
\begin{prop}\label{rig change} Let $P\in C_L^{(t)}$.

If $M^+_jP$ is defined if and only if $\la_j(P)<\la_{j-1}(P)$ and
in such case $\la_i(M^+_jP)=\la_i(P)+\delta_{i,j}$.

If $M^-_jP$ is defined if $\la_j(P)>\la_{j+1}(P)$ and
in such case  $\la_i(M^-_jP)=\la_i(P)-\delta_{i,j}$.
\end{prop}

\medskip
Now we describe paths with zero rigging. Fix an integer $m\geq1$.
The following lemma follows from the definitions.
\begin{lem}\label{0 rig}
There is a bijection from the set of rigged paths in $C_{L,r}^{(t)}$
which have $m$ particles and zero rigging $\la_1=\dots=\la_m=0$
to the set of rigged paths in $C_{L-2m,r}^{(t)}$ which have no  
particles.
The bijection maps the path
\be
P_{(m)}=\left(\begin{matrix}
r_L&\dots &r_{2m+1}&1      & 2&\dots   &1   & 2 &1    &2\ \ 1\\
   &\dots &\si_{2m+1}&\si_{2m} &0 &\dots   &v^{(t)}(1)& 0 & v^{(t)} 
(1)&0\ \ 0
\end{matrix}\right)\in C_{L,r}^{(t)}
\en
to the path
\be
P_{(0)}=\left(\begin{matrix}
r_L&\dots &r_{2m+1} & 1\\
   &\dots &\si_{2m+1} & \si_{2m}-v^{(t)}(1)
\end{matrix}\right)\in C_{L-2m,r}^{(t)}.
\en
\end{lem}

{\bf Examples.} We set $v=v^{(t)}(1)$.

$L=2,m(P)=1$:
\be
\begin{pmatrix}
1&2&1\\&0&0
\end{pmatrix}
\buildrel{M^+_1}\over\longrightarrow
\begin{pmatrix}
1&2&1\\&1&0
\end{pmatrix}
\buildrel{M^+_1}\over\longrightarrow
\begin{pmatrix}
1&2&1\\&0&1
\end{pmatrix}
\buildrel{M^+_1}\over\longrightarrow
\begin{pmatrix}
1&2&1\\&1&1
\end{pmatrix}
\buildrel{M^+_1}\over\longrightarrow\cdots
\en

$L=3,m(P)=1,a\geq0$:
\be
&&\begin{pmatrix}
2&1&2&1\\&v+a&0&0
\end{pmatrix}
\buildrel{M^+_1}\over\longrightarrow
\begin{pmatrix}
2&1&2&1\\&v+a-1&1&0
\end{pmatrix}
\buildrel{M^+_1}\over\longrightarrow
\begin{pmatrix}
2&1&2&1\\&v+a-1&0&1
\end{pmatrix}\\[2pt]
&\buildrel{M^+_1}\over\longrightarrow&
\begin{pmatrix}
2&1&2&1\\&v&0&a
\end{pmatrix}
\buildrel{M^+_1}\over\longrightarrow
\begin{pmatrix}
2&1&2&1\\&v-1&1&a
\end{pmatrix}
\buildrel{M^+_1}\over\longrightarrow
\cdots
\buildrel{M^+_1}\over\longrightarrow
\begin{pmatrix}
2&1&2&1\\&0&v&a
\end{pmatrix}\\[2pt]
&\buildrel{M^+_1}\over\longrightarrow&
\begin{pmatrix}
2&3&2&1\\&0&0&a
\end{pmatrix}
\buildrel{M^+_1}\over\longrightarrow
\begin{pmatrix}
2&1&2&1\\&0&v+1&a
\end{pmatrix}
\buildrel{M^+_1}\over\longrightarrow
\begin{pmatrix}
2&3&2&1\\&0&1&a
\end{pmatrix}
\buildrel{M^+_1}\over\longrightarrow\cdots
\en

$L=4,m(P)=1,a\geq0,2\leq b\leq v$:
\be
&&\begin{pmatrix}
1&2&1&2&1\\&b&v+a&0&0
\end{pmatrix}
\buildrel{M^+_1}\over\longrightarrow
\cdots
\buildrel{M^+_1}\over\longrightarrow
\begin{pmatrix}
1&2&1&2&1\\&b&v&0&a
\end{pmatrix}
\buildrel{M^+_1}\over\longrightarrow
\cdots
\buildrel{M^+_1}\over\longrightarrow
\begin{pmatrix}
1&2&1&2&1\\&b&v-b&b&a
\end{pmatrix}\\[2pt]
&&\buildrel{M^+_1}\over\longrightarrow
\cdots
\buildrel{M^+_1}\over\longrightarrow
\begin{pmatrix}
1&2&1&2&1\\&0&v&b&a
\end{pmatrix}
\buildrel{M^+_1}\over\longrightarrow
\begin{pmatrix}
1&2&1&2&1\\&1&v&b&a
\end{pmatrix}
\buildrel{M^+_1}\over\longrightarrow
\begin{pmatrix}
1&2&1&2&1\\&0&v+1&b&a
\end{pmatrix}
\buildrel{M^+_1}\over\longrightarrow
\cdots
\en

$L=4,m(P)=1,a\geq0,b>v$:
\be
&&\begin{pmatrix}
1&2&1&2&1\\&b&v+a&0&0
\end{pmatrix}
\buildrel{M^+_1}\over\longrightarrow
\cdots
\buildrel{M^+_1}\over\longrightarrow
\begin{pmatrix}
1&2&1&2&1\\&b&v&0&a
\end{pmatrix}
\buildrel{M^+_1}\over\longrightarrow
\cdots
\buildrel{M^+_1}\over\longrightarrow
\begin{pmatrix}
1&2&1&2&1\\&b&0&v&a
\end{pmatrix}\\[2pt]
&&\buildrel{M^+_1}\over\longrightarrow
\begin{pmatrix}
1&2&3&2&1\\&b-v-1&0&0&a
\end{pmatrix}
\buildrel{M^+_1}\over\longrightarrow
\begin{pmatrix}
1&2&1&2&1\\&b-1&0&v+1&a
\end{pmatrix}
\buildrel{M^+_1}\over\longrightarrow
\cdots\\[2pt]
&&\buildrel{M^+_1}\over\longrightarrow
\begin{pmatrix}
1&2&1&2&1\\&v&0&b&a
\end{pmatrix}
\buildrel{M^+_1}\over\longrightarrow
\cdots
\en

$L=4,m(P)=1,a,b\geq0$: Set $v'=v^{(t)}(2)$.
\be
&&\begin{pmatrix}
3&2&1&2&1\\&b&v+a&0&0
\end{pmatrix}
\buildrel{M^+_1}\over\longrightarrow
\cdots
\buildrel{M^+_1}\over\longrightarrow
\begin{pmatrix}
3&2&1&2&1\\&b&0&v&a
\end{pmatrix}
\buildrel{M^+_1}\over\longrightarrow
\begin{pmatrix}
3&2&3&2&1\\&v'+b&0&0&a
\end{pmatrix}\\[2pt]
&&\buildrel{M^+_1}\over\longrightarrow
\begin{pmatrix}
3&2&1&2&1\\&b-1&0&v+1&a
\end{pmatrix}
\buildrel{M^+_1}\over\longrightarrow
\cdots
\buildrel{M^+_1}\over\longrightarrow
\begin{pmatrix}
3&2&1&2&1\\&0&0&v+b&a
\end{pmatrix}\\[2pt]
&&\buildrel{M^+_1}\over\longrightarrow
\begin{pmatrix}
3&2&3&2&1\\&v'&0&b&a
\end{pmatrix}
\buildrel{M^+_1}\over\longrightarrow
\cdots
\buildrel{M^+_1}\over\longrightarrow
\begin{pmatrix}
3&2&3&2&1\\&0&v'&b&a
\end{pmatrix}
\buildrel{M^+_1}\over\longrightarrow
\begin{pmatrix}
3&4&3&2&1\\&0&0&b&a
\end{pmatrix}\\[2pt]
&&\buildrel{M^+_1}\over\longrightarrow
\begin{pmatrix}
3&2&3&2&1\\&0&v'+1&b&a
\end{pmatrix}
\buildrel{M^+_1}\over\longrightarrow
\begin{pmatrix}
3&4&3&2&1\\&0&1&b&a
\end{pmatrix}
\buildrel{M^+_1}\over\longrightarrow
\cdots
\en

$L=4,m(P)=2$:
\be
&&\scriptstyle{
\begin{pmatrix}
   \sr{1}&\sr{2}&\sr{1}&\sr{2}&\sr{1}\\
         &\sr{0}&\sr{v}&\sr{0}&\sr{0}
  \end{pmatrix}
\buildrel{M^+_1}\over\longrightarrow
  \begin{pmatrix}
   \sr{1}&\sr{2}&\sr{1}&\sr{2}&\sr{1}\\
         &\sr{1}&\sr{v}&\sr{0}&\sr{0}
  \end{pmatrix}
\buildrel{M^+_1}\over\longrightarrow
   \begin{pmatrix}
    \sr{1}&\sr{2}&\sr{1}&\sr{2}&\sr{1}\\
        &\sr{0}&\sr{v+1}&\sr{0}&\sr{0}
   \end{pmatrix}
\buildrel{M^+_1}\over\longrightarrow
   \begin{pmatrix}
    \sr{1}&\sr{2}&\sr{1}&\sr{2}&\sr{1}\\
        &\sr{1}&\sr{v+1}&\sr{0}&\sr{0}
   \end{pmatrix}
\buildrel{M^+_1}\over\longrightarrow\cdots}\\
&&\hskip125pt\sr{\downarrow M^+_2}
\hskip75pt\sr{\downarrow M^+_2}
\hskip75pt\sr{\downarrow M^+_2}\\
&&\scriptstyle{
\phantom{\begin{pmatrix}
   \sr{1}&\sr{2}&&\sr{2}&\sr{1}\\
         &\sr{0}&&\sr{0}&\sr{0}
  \end{pmatrix}
\buildrel{M^+_1}\over\longrightarrow}
  \begin{pmatrix}
   \sr{1}&\sr{2}&\sr{1}&\sr{2}&\sr{1}\\
         &\sr{1}&\sr{v-1}&\sr{1}&\sr{0}
  \end{pmatrix}
\buildrel{M^+_1}\over\longrightarrow\
   \begin{pmatrix}
    \sr{1}&\sr{2}&\sr{1}&\sr{2}&\sr{1}\\
        &\sr{0}&\sr{v}&\sr{1}&\sr{0}
   \end{pmatrix}
\buildrel{M^+_1}\over\longrightarrow\
   \begin{pmatrix}
    \sr{1}&\sr{2}&\sr{1}&\sr{2}&\sr{1}\\
        &\sr{1}&\sr{v}&\sr{1}&\sr{0}
   \end{pmatrix}
\ \buildrel{M^+_1}\over\longrightarrow\cdots}\\
&&\phantom{\hskip125pt\sr{\downarrow M^+_2}}
\hskip75pt\sr{\downarrow M^+_2}
\hskip75pt\sr{\downarrow M^+_2}\\
&&\scriptstyle{
\phantom{\begin{pmatrix}
   \sr{1}&\sr{2}&&\sr{2}&\sr{1}\\
         &\sr{0}&&\sr{0}&\sr{0}
  \end{pmatrix}
\buildrel{M^+_1}\over\longrightarrow
  \begin{pmatrix}
   \sr{1}&\sr{2}&\sr{1}&\sr{2}&\sr{1}\\
         &\sr{1}&\sr{v-1}&\sr{1}&\sr{0}
  \end{pmatrix}
\buildrel{M^+_1}\over\longrightarrow\ }
   \begin{pmatrix}
    \sr{1}&\sr{2}&\sr{1}&\sr{2}&\sr{1}\\
        &\sr{0}&\sr{v}&\sr{0}&\sr{1}
   \end{pmatrix}
\buildrel{M^+_1}\over\longrightarrow\
   \begin{pmatrix}
    \sr{1}&\sr{2}&\sr{1}&\sr{2}&\sr{1}\\
        &\sr{1}&\sr{v}&\sr{0}&\sr{1}
   \end{pmatrix}
\ \buildrel{M^+_1}\over\longrightarrow\cdots}\\
&&\phantom{\hskip125pt\sr{\downarrow M^+_2}
\hskip75pt\sr{\downarrow M^+_2}}
\hskip75pt\sr{\downarrow M^+_2}\\
&&\scriptstyle{
\phantom{\begin{pmatrix}
   \sr{1}&\sr{2}&&\sr{2}&\sr{1}\\
         &\sr{0}&&\sr{0}&\sr{0}
  \end{pmatrix}
\buildrel{M^+_1}\over\longrightarrow
  \begin{pmatrix}
   \sr{1}&\sr{2}&\sr{1}&\sr{2}&\sr{1}\\
         &\sr{1}&\sr{v-1}&\sr{1}&\sr{0}
  \end{pmatrix}
\buildrel{M^+_1}\over\longrightarrow\
   \begin{pmatrix}
    \sr{1}&\sr{2}&\sr{1}&\sr{2}&\sr{1}\\
        &\sr{0}&\sr{v}&\sr{0}&\sr{1}
   \end{pmatrix}
\buildrel{M^+_1}\over\longrightarrow\ }
   \begin{pmatrix}
    \sr{1}&\sr{2}&\sr{1}&\sr{2}&\sr{1}\\
        &\sr{1}&\sr{v-1}&\sr{1}&\sr{1}
   \end{pmatrix}
\buildrel{M^+_1}\over\longrightarrow\cdots}\\
\en

\subsection{Bijection}\label{subsec:Bijection}
We are in a position to construct the bijection which leads to  
Proposition
\ref{eq:recursion-str-2}.

Let $\bar P\in C_{L}^{(t-1)}$
\be
\bar P=\left(\begin{matrix}
\bar r_{L}& \bar r_{L-1}&\dots&\bar r_0\\
& \bar \si_{L-1}&\dots&\bar \si_0
\end{matrix}\right).\notag
\en

Define the numbers $\si_0,\dots,\si_{L-1}$ and the path $P_{(0)}\in  
C^{(t)}_L$
by the formulas
\be
\si_0&=&\bar \si_0,\\
\si_i&=&\bar \si_i +1/2+w^{(t-1)}(\bar r_{i+1},\bar r_i,\bar r_{i-1})-
w^{(t)}(\bar r_{i+1},\bar r_i,\bar r_{i-1}),
\quad(1\leq i\leq L-1),\\
P_{(0)}
&=&\left(\begin{matrix}
\bar r_{L}& \bar r_{L-1}&\dots&\bar r_0\\
&  \si_{L-1}&\dots& \si_0
\end{matrix}\right).
\en

\begin{lem}\label{IOTA0}
We have an inclusion
\be
\iota_0:\ C_{L,r}^{(t-1)} &\to &  C_{L,r}^{(t)}, \notag \\
\bar P &\mapsto & P_{(0)}.
\en
The image of $\iota_0$ coincides with the subset of paths in
$C_{L,r}^{(t)}$ with no particles.
\end{lem}
\begin{proof}
For $i>1$ by a simple computation we have

$\si_i=\bar \si_i$ if $r_{i-1}-r_{i+1}=\pm 2$,

$\si_i=\bar \si_i+1$ if $r_{i-1}-r_{i+1}=0$ and $r_{i-1}\neq 1,p-1$,

$\si_i=\bar \si_i+2$ if $r_{i-1}=r_{i+1}=1$ or  $r_{i-1}=r_{i+1}=p-1$.

\noindent The statement of the lemma follows straightforwardly.
\end{proof}

Let $\bar P\in C_{L-2m,r}^{(t-1)}$,
and let $\la=(\la_1,\dots,\la_m)\in \pi_m$
be a partition. The path $P_{(0)}=\iota_0(\bar P)\in C_{L-2m,r}^{(t)}$
is constructed in Lemma \ref{IOTA0}. The path $P_{(m)}\in C_{L,r}^ 
{(t)}$
is constructed from $P_{(0)}$ in Lemma \ref{0 rig}.
\begin{lem}\label{m incl}
Notation being as above,
for any non-negative integer $m$ we have an inclusion
\be
\iota_m:\ C_{L-2m,r}^{(t-1)}\times \pi_m &\to &  C_{L,r}^{(t)},  
\notag \\
(\bar P,\lambda) &\mapsto &
(M^+_m)^{\la_m} \dots (M^+_2)^{\la_2}  (M^+_1)^{\la_1}P_{(m)}.
\en
The image of $\iota_m$ coincides with the subset of paths in
$C_{L,r}^{(t)}$ with $m$ particles. We have
\bea\label{DEG}
d(\iota_m(\bar P,\la))=d(\bar P)+\sum_{i=1}^m \la_i+L^2/4+L/2\,.
\ena
\end{lem}
\begin{proof}
The degree of a path is defined in (\ref{DEGREE}).
The relation (\ref{DEG}) follows by a straightforward computation.
The rest follows from Lemma \ref{rig change}.
\end{proof}

We obtain the main result of this section.
\begin{thm}\label{thm:iota-is-bijection}
The map $\iota:=\sqcup_{m=0}^\infty \iota_m$ defines a bijection
\be
\bigsqcup_{m=0}^\infty (C_{L-2m,r}^{(t-1)}\times \pi_m) \to C_{L,r}^ 
{(t)}
\en
with the property $(\ref{DEG})$.
\end{thm}
Proposition \ref{prop:ch_t_to_t-1} is a direct consequence of
Theorem \ref{thm:iota-is-bijection}.
\setcounter{section}{4}
\setcounter{equation}{0}
\newcommand{\La}{\Lambda}
\newcommand{\End}{\mathop{\rm End}}
\newcommand{\at}{\tilde{a}}
\newcommand{\ft}{\tilde{f}}
\newcommand{\Vt}{\widetilde{V}}
\newcommand{\Wt}{\widetilde{W}}
\newcommand{\vt}{{\tilde{v}}}
\newcommand{\phit}{\tilde{\phi}}
\newcommand{\id}{{\rm id}}
\newcommand{\ut}{{\tilde{u}}}
\newcommand{\Sn}{\mathfrak{S}}
\section{The case $p=3$}\label{sec:p=3}
The aim of this section is to deduce
Conjecture \ref{CONJ} for $p=3$
{}from the work \cite{FJM}.

Throughout this section, we fix an integer $p'>3$ coprime to $3$.
Omitting the upper index we write the $(2,1)$-field
$\phi^{(r',r)}(z)$ as $\phi(z)$, since for $p=3$ the choice of $r'=r 
\pm1$
is uniquely determined from $r$. We set
\be
M^{(3,p')}=M^{(3,p')}_{1,1}\oplus M^{(3,p')}_{2,1}.
\en
In the following we deal with
bi-graded $\C$-vector spaces of the form
$X=\oplus_{d,L}X_{d,L}$ with $\dim X_{d,L}<\infty$.
We will refer to the index $d$ and $L$ as
degree and weight, respectively.
We set $X_{L}=\oplus_dX_{d,L}$ and define the restricted dual space by
$X^*_L=\oplus_d\Hom_\C(X_{d,L},\C)$.

\subsection{Extended modules and monomial basis}
\label{subsec:mon1}
First we review the results of \cite{FJM}
which are relevant to us.

Consider the Heisenberg algebra with generators
$\{h_n\}_{n\in\Z}$ satisfying $[h_m,h_n]=m\delta_{m+n,0}$.
Denote by
\be
\F_{\gamma}=\C[h_{-1},h_{-2},\cdots]\ket{\gamma}
\en
the Fock space with highest weight vector
$\ket{\gamma}$ ($\gamma\in\C$),
where $h_n\ket{\gamma}=0$ ($n>0$)
and $h_0\ket{\gamma}=\gamma\ket{\gamma}$.
We set
\be
\beta=\sqrt{\frac{p'-2}{2}}.
\en
We use the vertex operator which acts on
$\F=\oplus_{n\in\Z}\F_{n\beta}$,
\be
\Phi_\beta(z)=
:\exp\left(-\beta\sum_{n\neq0}\frac{h_n}{n}z^{-n}\right)
e^{\beta Q}z^{\beta h_0}:.
\en
Here $:\phantom{\phi}:$ stands for the normal ordering symbol,
and $e^{\beta Q}:\F_\gamma\overset{\sim}{\to}\F_{\gamma+\beta}$
is an isomorphism of vector spaces such that
$[h_n,e^{\beta Q}]=0$ ($n\neq 0$),
$e^{\beta Q}\ket{\gamma}=\ket{\gamma+\beta}$.
Set
\be
V^{(3,p')}&=&\oplus_{L\in\Z}V^{(3,p')}_L\hbox{ where }
V^{(3,p')}_L=
\begin{cases}
M^{(3,p')}_{1,1}\otimes\F_{L\beta}&\hbox{ if $L$ is even};\\
M^{(3,p')}_{2,1}\otimes\F_{L\beta}&\hbox{ if $L$ is odd},
\end{cases}\\
\vac
&=&\ket{1,1}\otimes\ket{0}~~\in V^{(3,p')}_0.
\en
Introduce a field with coefficients in $\End(V^{(3,p')})$,
\be
a(z)=\phi(z)\otimes \Phi_{\beta}(z).
\en
The parameter $\beta$ is so chosen that $a(z)$ has the expansion
\be
a(z)=\sum_{n\in\Z}a_{n}z^{-n-1},
\en
and that the coefficients $a_n$'s are mutually commutative,
\bea
[a_m,a_n]=0 \qquad (m,n\in\Z).
\label{eq:a-comm}
\ena
We have also $a_n\vac=0$ for $n\ge 0$.
Later we will use the relation
\bea
a(z)^2=k \cdot \id \, \otimes :\Phi_\beta(z)^2:,
\label{eq:a2=0}
\ena
where $k$ is a nonzero constant.
We have
\bea
\label{T1}
[h_n,a(z)]=\beta z^na(z).
\ena

Let $A$ be a subalgebra of ${\rm End}(V^{(3,p')})$ generated by
$a_n$'s.
Our main concern is the following subspace of $V^{(3,p')}$  
generated from
$\vac$ by acting with $a_n$'s:
\be
W^{(3,p')}=A\cdot\vac\,
\subset V^{(3,p')}\,.
\en
We introduce a bi-grading to $W^{(3,p')}$ by
assigning the degree and weight as
$\deg a_n=-n$, $\wt a_n=1$
and $\deg \vac=0$, $\wt \vac=0$.

Let $\La_L=\C[x_1,\cdots,x_L]^{\Sf_L}$
denote the space of symmetric
polynomials in the variables $x_1,\cdots,x_L$.
The restricted dual space $W^{(3,p')*}_L$ is identified with a  
subspace of
$\Lambda_L$ by
\bea
W^{(3,p')*}_L\longrightarrow \Lambda_L,
\quad
\langle v| \mapsto \langle v|a(x_L)\cdots a(x_1)\vac\,.
\label{eq:dualpol}
\ena
If $n>0$, $h_n\vac=0$ and from \eqref{T1} we see that
$W^{(3,p')}_L$ is invariant by $h_n$.
The right action of $h_n$ ($n>0$) on $W^{(3,p')*}_L$
corresponds to the multiplication by $\beta\sum_{j=1}^L x_j^n$ on $ 
\La_L$.
Therefore,
the image $I^{(3,p')}_L$ of \eqref{eq:dualpol} is an ideal of $\La_L$.

The linear isomorphism
$\tau=\id\otimes e^{2\beta Q}:V^{(3,p')}\to V^{(3,p')}$
has an effect of shifting the indices
\bea
\tau \vac_{2l}=\vac_{2l+2},\quad
\tau\,a_n\,\tau^{-1}=a_{n-p'+2},
\label{eq:shift}
\ena
where $\vac_{2l}=\ket{1,1}\otimes\ket{2l\beta}$.
We have an increasing filtration (\cite{FJM}, Proposition 2.3)
\bea
&&W^{(3,p')}\subset \tau^{-1}(W^{(3,p')})\subset
\tau^{-2}(W^{(3,p')})\subset  \cdots,
\label{eq:increase}\\
&&V^{(3,p')}=\bigcup_{N\ge 0}\tau^{-N}(W^{(3,p')}).
\label{eq:Wfilt}
\ena
Let $\widehat{\La}_L=\C[x_1^{\pm1},\cdots,x_L^{\pm1}]^{\Sf_L}$
be the space of symmetric Laurent
polynomials, and let $\widehat{I}^{(3,p')}_L$ be the ideal of
$\widehat{\La}_L$ generated by $I^{(3,p')}_L$.
By \eqref{eq:a-comm}, \eqref{eq:shift} and \eqref{eq:Wfilt}, we see  
that
the subspace of $\widehat{\La}_L$ spanned by all matrix elements
\bea
\bra{v}a(x_L)\cdots a(x_1)\ket{u}
\quad (\bra{v}\in (V^{(3,p')})^*, \ket{u}\in V^{(3,p')})
\label{eq:mat-elt}
\ena
coincides with $\widehat{I}^{(3,p')}_L$.
Note also that $f\in\widehat{I}^{(3,p')}_L$ if and only if
$x_1\cdots x_Lf\in\widehat{I}^{(3,p')}_L$.

The following fact is proved in \cite{FJM}, Theorem 2.8.
\begin{prop}\label{prop:principal}
For each $L\ge 0$, the space
$W^{(3,p')}_L$ has a basis consisting of elements
\be
a_{-\la_L}\cdots a_{-\la_1}\vac,
\en
where
$\la=(\la_L,\cdots,\la_1)$ runs over all $L$-tuples of positive
integers satisfying
\bea
&&\la_{L}\ge\cdots\ge \la_1,
\label{eq:adm1}\\
&&\la_{i+2}-\la_{i}\ge p'-2
\quad (1\le i\le L-2).
\label{eq:mon-a2}
\ena
\end{prop}

Let us reformulate Proposition \ref{prop:principal}
in a form suitable for later use.
\begin{prop}\label{prop:cubic-rel}
For any $\la=(\la_L,\cdots,\la_1)\in\Z^L$
with $\la_L\ge\cdots\ge\la_1$,
there exist unique $c_{\la\mu}\in\C$ such that
the following identity holds as operators on $V^{(3,p')}:$
\bea
a_{-\la_L}\cdots a_{-\la_1}
=\sum_{\mu}c_{\la\mu}\,
a_{-\mu_L}\cdots a_{-\mu_1},
\label{eq:operator-id}
\ena
where the sum in the right hand side
is taken over all $\mu=(\mu_L,\cdots,\mu_1)\in\Z^L$
satisfying \eqref{eq:adm1}, \eqref{eq:mon-a2} and
$\sum_{i=1}^L\mu_i=\sum_{i=1}^L\la_i$. The coefficients $c_{\la,\mu}$
are shift invariant\,$:$
\be
c_{\la+(1,\ldots,1),\mu+(1,\ldots,1)}=c_{\la,\mu}.
\en
\end{prop}
\begin{proof}
The sum in the right hand side of \eqref{eq:operator-id}
is a well defined element of ${\rm End}(V^{(3,p')})$
because for a given vector in $V^{(3,p')}$ only finitely many
monomials $a_{-\mu_L}\cdots a_{-\mu_1}$ acts non-trivially.
For $N\in\Z$ and $\mu=(\mu_L,\cdots,\mu_1)\in\Z^L$, let us write
$\mu^{(N)}=(\mu^{(N)}_L,\cdots,\mu^{(N)}_1)$, $\mu_i^{(N)}=\mu_i+N 
(p'-2)$.
Choosing $N$ so that $\la^{(N)}_1>0$, we apply
Proposition \ref{prop:principal} to $\la^{(N)}$.
There exist unique constants $c^{(N)}_{\la\mu}\in\C$ such that
\be
\Bigl(a_{-\la^{(N)}_L}\cdots a_{-\la^{(N)}_1}-
\sum_{\mu}c^{(N)}_{\la\mu}\,
a_{-\mu^{(N)}_L}\cdots a_{-\mu^{(N)}_1}\Bigr)
\vac=0,
\en
where the sum is taken over  $\mu$
satisfying \eqref{eq:adm1}, \eqref{eq:mon-a2}
and $\mu_1>-N(p'-2)$.
Since $a_n$'s commute and $W^{(3,p')}$ is cyclic over
$\C[a_{-1},a_{-2},\cdots]$,
the operator inside the parentheses vanishes on  $W^{(3,p')}$.
Applying $\tau^{-N}$ on both sides, we deduce that
\be
a_{-\la_L}\cdots a_{-\la_1}
=\sum_{\mu}c^{(N)}_{\la\mu}a_{-\mu_L}\cdots a_{-\mu_1}
\quad \mbox{on $\tau^{-N}(W^{(3,p')})$}.
\en
If $N'>N$, then $\tau^{-N}(W^{(3,p')})\subset \tau^{-N'}(W^{(3,p')})$.
Therefore we have
\be
\sum_\mu(c^{(N)}_{\la\mu}-c^{(N')}_{\la\mu})
a_{-\mu_L}\cdots a_{-\mu_1}=0
\quad \mbox{on $\tau^{-N}(W^{(3,p')})$},
\en
or equivalently
\be
\sum_\mu(c^{(N)}_{\la\mu}-c^{(N')}_{\la\mu})
a_{-\mu_L-N(p'-2)}\cdots a_{-\mu_1-N(p'-2)}=0
\quad \mbox{on $W^{(3,p')}$}.
\en
{}From Proposition \ref{prop:principal}
we have $c^{(N')}_{\la\mu}=c^{(N)}_{\la\mu}$ for all
$\mu$ with $\mu_1> -N(p'-2)$.
Setting $c_{\la\mu}=\lim_{N\to\infty}c^{(N)}_{\la\mu}$
and noting \eqref{eq:Wfilt}, we obtain \eqref{eq:operator-id}.
The uniqueness of $c_{\la\mu}$ is clear.

The last statement follows from the commutation relation \eqref{T1}.
\end{proof}
\begin{cor}\label{cor:f=0}
A symmetric Laurent polynomial
$f=\sum_{\la_1,\cdots,\la_L\in\Z}f_{\la_L,\cdots,\la_1}
x_L^{\la_L}\cdots x_1^{\la_1}$
belongs to $\widehat{I}^{(3,p')}_L$ if and only if
\be
f_{\la_L,\ldots,\la_1}=
\sum_{\mu}c_{\la\mu}f_{\mu_L,\ldots,\mu_1}
\en
holds for all $\la=(\la_L,\cdots,\la_1)\in\Z^L$ with $\la_L\ge\cdots 
\ge\la_1$.
In particular, if $f_{\la_L,\cdots,\la_1}=0$ holds
for all $\la$ satisfying \eqref{eq:adm1}, \eqref{eq:mon-a2}, then  
$f=0$.
\end{cor}

We quote one more fact which
follows immediately from \cite{FJM}, Propositions 4.2 and 4.4.
\begin{prop}\label{prop:dual-pol}
Let $\widehat{D}_L$ be the subspace of $\widehat{\Lambda}_L$
consisting of Laurent polynomials
vanishing on the diagonal $x_i-x_j=0$ $(i\neq j)$.
Then
\be
\widehat{I}^{(3,p')}_L\cap \widehat{D}_L=
\widehat{I}^{(3,p'-3)}_L\times\prod_{1\le i<j\le L}(x_i-x_j)^2.
\en
When $p'=4,5$, we define $\widehat{I}^{(3,p'-3)}_L$
in the right hand side to be $\widehat{\La}_L$.
\end{prop}


\subsection{Monomial basis for $M^{(3,p')}$}\label{subsec:mon-2}
For $p=3$, the admissibility of the rigged paths is equivalent to
the following conditions for the $n_i$'s (see Section 2.3).
\bea
&&n_1\in\Z_{\geq0}+\Delta_{2,1},\label{n1}\\
&&n_{i+1}-n_i\ge -\frac{p'}{2}+3\quad (1\le i\le L-1),
\label{eq:adm-3p1}\\
&&n_{i+2}-n_i\ge 1\quad (1\le i\le L-2).
\label{eq:adm-3p2}
\ena
Our goal in this section is to prove the following result.
\begin{thm}\label{thm:p=3}
The set of monomials
\be
\phi_{-n_L}\cdots\phi_{-n_1}\ket{1,1}
\en
satisfying the conditions
\eqref{n1}, \eqref{eq:adm-3p1} and \eqref{eq:adm-3p2}
is a basis of $M^{(3,p')}$.
\end{thm}
We have shown in Theorem \ref{thm:main-theorem} that
the character of this set matches with that of $M^{(3,p')}$.
For the proof of Theorem \ref{thm:p=3}, it is therefore sufficient
to show that this set spans $M^{(3,p')}$.
For that purpose, we consider a filtration of $M^{(3,p')}$
\be
\{0\}=F_{-1}\subset F_0\subset\cdots\subset F_i\subset
\cdots\subset M^{(3,p')},
\en
defined by
\bea
&&F_0=\C\ket{1,1},
\label{eq:F-filt1}\\
&&F_i=F_{i-1}+\sum_{n\in\Z
+(-1)^{i-1}\Delta_{2,1}}\phi_{-n}F_{i-1},
\label{eq:F-filt2}
\ena
where $\Delta_{2,1}=(p'-2)/4$.
In the right hand side of \eqref{eq:F-filt2}, we have
\be
\phi_{-n}=
\begin{cases}
\phi^{(2,1)}_{-n}&\hbox{if $i$ is odd};\\
\phi^{(1,2)}_{-n}&\hbox{if $i$ is even}.\\
\end{cases}
\en
If $r\neq r''$
the action of $\phi_{-n}^{(r',r)}$ is zero
on $M^{(3,p')}_{r'',1}$ by definition.
Since the operator product expansion of $\phi(z)$ contains the
energy-momentum tensor $T(z)$, we have $\cup_{i\geq0}F_i=M^{(3,p')}$.

Let
\be
\phit_n:F_{i}/F_{i-1}\longrightarrow F_{i+1}/F_i
\en
denote the operator induced from $\phi_n$ on the associated graded  
space $\widetilde M^{(3,p')}=\oplus_{i\ge 0}F_{i}/F_{i-1}$.
{}From the remark made above, the proof of
Theorem \ref{thm:p=3} reduces to the following statement.
\begin{lem}\label{lem:span-Wt}
Let $\widetilde{B}^{(3,p')}_L$ denote the set of vectors
\bea
\phit_{-n_L}\cdots\phit_{-n_1}\ket{1,1}
\qquad (n_i\in\Z-(-1)^i\Delta_{2,1})
\label{eq:mono-Wt}
\ena
satisfying the conditions
\eqref{n1}, \eqref{eq:adm-3p1} and \eqref{eq:adm-3p2}.
Then $\widetilde M^{(3,p')}_L={\rm span}_\C\widetilde{B}^{(3,p')}_L$.
\end{lem}
\noindent
The rest of this subsection is devoted to the proof of
Lemma \ref{lem:span-Wt}.

We set
\be
\phit(z)=\sum_{n\in\Z+(-1)^i\Delta_{2,1}}
\phit_{-n} z^{n-\Delta_{2,1}}
\quad \mbox{on $F_{i}/F_{i-1}$}.
\en
The field $\phit(z)$
satisfies a quadratic relation following from that of $\phi(z)$
(see e.g.\cite{FJMMT}),
\be
\sum_{k\ge0}c_k\,\phit_{-n_2-k}\phit_{-n_1+k}
+
\sum_{k\ge0}c_k\,\phit_{-n_1+p'/2-2-k}\phit_{-n_2-p'/2+2+k}=0,
\en
where $\sum_{k\ge 0}c_kz^k=(1-z)^{p'/2-2}$.
This can be rewritten as a relation of the form
\bea
\phit_{-n_2}\phit_{-n_1}
=\sum_{l\ge 1}C_{n_1,n_2,l}\phit_{-n_2-l}\phit_{-n_1+l}
\quad (n_2-n_1\le -\frac{p'}{2}+2)
\label{eq:quad-rel}
\ena
with some $C_{n_1,n_2,l}\in\C$.
Given any monomial \eqref{eq:mono-Wt}, we apply
\eqref{eq:quad-rel} to rewrite it as a linear combination of
monomials satisfying the condition \eqref{eq:adm-3p1}.
In each step of rewriting,
$\sum_{i=1}^Ln_i$ is invariant and
$\sum_{i=1}^L in_i$ strictly increases.
Since $\sum_{j=1}^i n_j>0$ for each $i$, $\sum_{i=1}^L in_i$
is bounded from above. Therefore, the process terminates after a
finite number of steps.
To prove Lemma \ref{lem:span-Wt}, we must also fulfill
the second condition \eqref{eq:adm-3p2}.
To this end, we make use of the knowledge in the previous subsection.

Consider the filtration of $V^{(3,p')}$ induced from
that of $M^{(3,p')}$.
\be
V^{(3,p')}_i&=&\oplus_{L\in\Z}(V^{(3,p')}_L)_i,\\
(V^{(3,p')}_L)_i&=&
\begin{cases}
(F_i\cap M^{(3,p')}_{1,1})\otimes \F_{L\beta}&\hbox{ if $L$ is  
even};\\
(F_i\cap M^{(3,p')}_{2,1})\otimes \F_{L\beta}&\hbox{ if $L$ is odd}.
\end{cases}
\en
Define
\be
\at_n:V^{(3,p')}_i/V^{(3,p')}_{i-1}\longrightarrow V^{(3,p')}_{i+1}/ 
V^{(3,p')}_i
\en
to be the operator induced from $a_n:V^{(3,p')}_i\to V^{(3,p')}_{i 
+1}$.
Setting
$\at(z)=\sum_{n\in\Z}\at_nz^{-n-1}$,
we have clearly the relation
\be
\at(z)=\phit(z)\otimes \Phi_\beta(z).
\en
Since the relations \eqref{eq:operator-id} are homogeneous,
they remain valid for $\tilde a_n$'s.
{}From the relation \eqref{eq:a2=0}, we have also $\at(z)^2=0$.

Now for $\bra{\vt}\in {\widetilde M^{(3,p')*}}$ and
$\ket{\ut}\in \widetilde M^{(3,p')}$,
we consider a matrix element
\be
&&
g(x_1,\ldots,x_L)=
\bra{\vt}\phit(x_L)\cdots\phit(x_1)\ket{\ut}
\en
which is a formal series in $x_1,\cdots,x_L$.

\begin{lem}\label{lem:g-f}
There exists an element $f\in \widehat{I}^{(3,p'-3)}_L$ such that
\be
g(x_1,\cdots,x_L)=f(x_1,\cdots,x_L)\times
\prod_{1\le i<j\le L}(x_j-x_i)^{-\beta^2+2}.
\en
Here the right hand side means its expansion in the domain
$|x_L|>\cdots>|x_1|$.
\end{lem}
\begin{proof}
Set $\bra{\vt_L}=\bra{\vt}\otimes\bra{L\beta}$,
$\ket{\ut_0}=\ket{\ut}\otimes\ket{0}$, and consider
\be
F(x_1,\cdots,x_L)
&=&\bra{\vt_L}\at(x_L)\cdots \at(x_1)\ket{\ut_0}
\\
&=&g(x_1,\cdots,x_L)\times
\prod_{1\le i<j\le L}(x_j-x_i)^{\beta^2}.
\en
We have $F\in \widehat{I}^{(3,p')}_L$ because the defining relations
\eqref{eq:operator-id} for the operator algebra dual to
the ideal $F\in \widehat{I}^{(3,p')}_L$ are equally valid for $ 
\tilde a(z)$.
Moreover the relation
$\at(z)^2=0$ implies $F\in\widehat{D}_L$.
Hence from Proposition \ref{prop:dual-pol}, we can write
\be
F(x_1,\cdots,x_L)
=f(x_1,\cdots,x_L)\prod_{1\le i<j\le L}(x_j-x_i)^{2}
\en
with some $f\in \widehat{I}^{(3,p'-3)}_L$.
This proves the assertion.
\end{proof}

\begin{lem}\label{lem:cubic-phi}
Let $\ket{\ut}\in \widetilde M^{(3,p')}$  be a
homogeneous element and let $m=(m_3,m_2,m_1)$ be a triple of integers
such that the relations \eqref{eq:adm-3p1}, \eqref{eq:adm-3p2} for  
$L=3$
are not satisfied. Then we have a relation
\bea
\phit_{-m_3}\phit_{-m_2}\phit_{-m_1}\ket{\ut}
=\sum_{n} d_{m,n}\,\phit_{-n_3}\phit_{-n_2}\phit_{-n_1}\ket{\ut}
\label{eq:cubic-phi}
\ena
with some $d_{m,n}\in\C$, where the sum is taken over
$n=(n_3,n_2,n_1)$ satisfying the relations
\eqref{eq:adm-3p1}, \eqref{eq:adm-3p2}
and $\sum_{i=1}^3n_i=\sum_{i=1}^3m_i$,
$\sum_{i=1}^3 i n_i>\sum_{i=1}^3 im_i$.
\end{lem}
\begin{proof}
Given $\ket{\ut}$, let
$X$ be the linear span of all monomials
$\phit_{-m_3}\phit_{-m_2}\phit_{-m_1}\ket{\ut}$,
and let $X_0$ be the subspace spanned by those
whose indices $m$ satisfy \eqref{eq:adm-3p1} and \eqref{eq:adm-3p2}.
We show $X=X_0$.

Let $\bra{\vt}\in {\widetilde M^{(3,p')*}}$
be an element orthogonal to the space $X_0$.
Consider the matrix element
$g(x_1,x_2,x_3)=
\bra{\vt}\phit(x_3)\phit(x_2)\phit(x_1)\ket{\ut}$.
Let $f$ be as in Lemma \ref{lem:g-f} with $L=3$, and expand them as
\be
&&
f(x_1,x_2,x_3)=\sum_{\la_1,\la_2,\la_3\in\Z}
f_{\lambda_3,\lambda_2,\lambda_1}x_3^{\lambda_3}
x_2^{\lambda_2}x_1^{\lambda_1},
\\
&&
g(x_1,x_2,x_3)=\sum_{n_1,n_2,n_3}
g_{n_3,n_2,n_1}x_3^{n_3-\Delta_{2,1}}x_2^{n_2-\Delta_{2,1}}x_1^{n_1- 
\Delta_{2,1}},
\en
where $g_{n_3,n_2,n_1}=
\langle\tilde v|\phit_{-n_3}\phit_{-n_2}\phit_{-n_1}\ket{\ut}$.

Suppose $\la_3\ge\la_2\ge\la_1$ and $\la_3-\la_1\ge p'-5$.
The coefficient $f_{\la_3,\la_2,\la_1}$ is a linear combination of
$g_{n_3,n_2,n_1}$ such that
\be
\la_i=n_i-\Delta_{2,1}+(i-1)(\beta^2-2)-\alpha_{i-1}+\alpha_i,
\en
with some $\alpha_1,\alpha_2\ge 0$ and $\alpha_0=\alpha_3=0$.
In particular
\be
&&n_3-n_1=\la_3-\la_1-(p'-6)+\alpha_1+\alpha_2\ge 1.
\en
Using the quadratic relation \eqref{eq:quad-rel},
we can further rewrite $g_{n_3,n_2,n_1}$ as linear combinations
of those satisfying \eqref{eq:adm-3p1}.
Since $n_3-n_1=\sum_{i=1}^3in_i-2\sum_{i=1}^3n_i$ does not decrease,
the resulting terms all satisfy \eqref{eq:adm-3p2} as well.
By the choice of $\bra{\vt}$,
$g_{n_3,n_2,n_1}=0$ holds for all such $n=(n_3,n_2,n_1)$.
It follows that $f_{\la_3,\la_2,\la_1}=0$
holds for all $\la=(\la_3,\la_2,\la_1)$
with  $\la_3\ge\la_2\ge\la_1$ and $\la_3-\la_1\ge p'-5$.
Applying Corollary \ref{cor:f=0} we conclude that $f=0=g$.
This shows $X=X_0$.

It remains to verify that if $m$ violates
\eqref{eq:adm-3p1} or \eqref{eq:adm-3p2}, then
$m_3-m_1<n_3-n_1$ for all terms
in the right hand side of \eqref{eq:cubic-phi}.
This is evident if \eqref{eq:adm-3p2} is violated.
Otherwise we can use the quadratic relation \eqref{eq:quad-rel} to  
rewrite
it as a linear combination of those satisfying \eqref{eq:adm-3p1}.
In the process $m_3-m_1$ strictly increases, and the verification  
reduces to
the first case.
\end{proof}
\medskip

\noindent{\it Proof of Lemma \ref{lem:span-Wt}. }\quad
Given a monomial \eqref{eq:mono-Wt},
suppose the conditions \eqref{eq:adm-3p1}--\eqref{eq:adm-3p2} are  
not valid
for a triple $(m_{i+1},m_i,m_{i-1})$.
We apply Lemma \ref{lem:cubic-phi} to reduce it.
In the process, $\sum_{i=1}^Ln_i$
does not change and $\sum_{i=1}^Lin_i$ strictly increases.
Therefore the process terminates after a finite number of steps.
Proof of Lemma \ref{lem:span-Wt} is now complete.
\qed

\bigskip
\noindent
{\it Acknowledgments.}\quad
BF is partially supported by grants RFBR-02-01-01015,
RFHR-01-01-00906, INTAS-00-00055.
MJ is partially supported by
the Grant-in-Aid for Scientific Research B2, no. 16340033,
and TM is partially supported by
(A1) no.13304010, Japan Society for the Promotion of Science.
EM is partially supported by the National Science
Foundation (NSF) grant DMS-0140460.
YT is partially supported by
Grant-in-Aid for Young Scientists (B) No.\ 17740089.

\end{document}